\documentclass[conference]{new-aiaa}

\usepackage[utf8]{inputenc}

\usepackage[version=4]{mhchem}
\usepackage{siunitx}

\usepackage{booktabs}
\usepackage{array}
\usepackage{longtable}
\usepackage{tabularx}

\usepackage{graphicx}
\usepackage{subcaption}
\usepackage{float}
\usepackage{rotating}

\usepackage{algorithm}
\usepackage{algpseudocode}

\usepackage{ifthen}
\usepackage{enumitem}
\usepackage{etoolbox}
\usepackage{xpatch}
\usepackage{xcolor}
\usepackage{soul}
\usepackage{textcomp}
\usepackage{chngcntr}
\usepackage{afterpage}
\usepackage{floatpag}

\usepackage{nomencl}
\makenomenclature

\title{Reliability, Robustness, and Resilience Modeling for Surveillance System in Advanced Air Mobility Operations}

\author{
Esrat Farhana Dulia\textsuperscript{1}\footnote{Graduate Research Assistant, College of Aeronautics and Engineering, edulia@kent.edu.}, 
Caleb Adams\textsuperscript{1}\footnote{Graduate Research Assistant, College of Aeronautics and Engineering, cadams68@kent.edu.}, 
Syed Arbab Mohd Shihab\textsuperscript{1}\footnote{Assistant Professor, College of Aeronautics and Engineering, sshihab@kent.edu.}, and
Ruben Del Rosario\textsuperscript{1}\footnote{Director, Center for Advanced Air Mobility and Professor, College of Aeronautics and Engineering, rdelros1@kent.edu.}
}

\affil{
\textsuperscript{1}Kent State University, Kent, OH, USA 44242 \\
}

\begin{document}

\maketitle

\begin{abstract}

Ensuring the safe and efficient operation of Advanced Air Mobility (AAM) in low-altitude airspace requires a reliable, robust, and resilient surveillance system capable of continuously detecting, identifying, and tracking aircraft under both normal and off-nominal conditions. To address this need, this study develops a comprehensive 3R modeling framework—reliability, robustness, and resilience—for the Surveillance for Advanced Air Mobility (SAM) system, with a focus on the optimal design and operation of a multi-type sensor network. Under normal operating conditions, the reliability model determines the baseline sensor types, quantities, and locations required to satisfy surveillance coverage and detection requirements. To address external perturbations, such as adverse weather conditions or sudden increases in AAM traffic demand, the robustness model identifies additional sensor requirements needed to maintain system performance. Furthermore, for surveillance outages caused by primary sensor failures, the resiliency model develops backup sensor deployment and dispatch strategies to provide temporary surveillance coverage, minimize operational disruptions, and support the safe continuation of AAM operations.

\end{abstract}

\section{Nomenclature}

\begin{longtable}{l l}
$\mathcal{I}$ & Set of candidate sensor sites \\ 
$\mathcal{S}$ & Set of sensor types \\ 
$\mathcal{K}$ & Set of AAM aircraft \\ 
$\mathcal{T}$ & Set of time steps in the planning horizon \\[3pt]
$\mathcal{H}$ & Set of candidate hub locations for backup sensors \\
$\mathcal{O}$ & Set of primary sensor locations \\
$\mathcal{B}$ & Set of backup sensors (UAVs or ground vehicles) \\
$n_{i}^{s}$ & Number of installed sensor sets of type $s$ at site $i$ \\ 
$\beta_{i}$ & Binary variable; 1 if site $i$ is activated, 0 otherwise \\ 

$c_{s}$ & Cost of installing one sensor of type $s$ \\ 
$C_{s}$ & Number of sensors contained in one set of type $s$ \\ 
$M$ & Maximum sensor capacity per site \\ 
$M_{s}$ & Maximum allowable number of sets of type $s$ at site $i$ \\[3pt]

$\lambda_{s}$ & Failure rate of sensor type $s$ \\ 
$\lambda_{\ell}$ & Failure rate of communication link \\ 
$\lambda_{u}$ & Failure rate of data server  \\ 

$\rho_{s}(t)$ & Reliability of sensor type $s$ at time $t$ \\ 
$\rho_{\ell}(t)$ & Reliability of communication link at time $t$ \\ 
$\rho_{u}(t)$ & Reliability of server at time $t$ \\[3pt]

$d_{i,k}(t)$ & 3D distance between site $i$ and AAM aircraft $k$ at time $t$ \\ 
$r_{s}$ & Effective sensing range of sensor type $s$ \\ 
$\chi_{i,s,k}(t)$ & Binary range indicator; 1 if $d_{i,k}(t)\le r_s$ \\ 
$\ell_{i,s,k}(t)$ & LOS quality factor between sensor $(i,s)$ and AAM aircraft $k$ at time $t$ \\[3pt]

$p_{i,s,k}(t)$ & Intrinsic detection probability for sensor $(i,s)$ and AAM aircraft $k$ at time $t$ \\ 
$q_{i,s,k}(t)$ & Probability that detection from sensor $(i,s)$ is successfully sent to data server \\ 
$m_{i,s,k}(t)$ & Missed detection probability for sensor $(i,s)$ \\ 
$\varepsilon$ & Small positive constant for numerical stability \\[3pt]

$\alpha_{t,i,s}$ & Set of AAM aircraft detected and reported by sensor type $s$ at site $i$ at time $t$ \\[3pt]

$H$ & Minimum required system detection reliability  \\[3pt]

$Z$ & Total cost of sensor deployment (objective function) \\ 
$n_{i,s}^{\text{exist}}$ & Number of existing sensors of type $s$ at site $i$ \\
$n_{i,s}^{\text{add}}$ & Number of additional sensors of type $s$ to be deployed at site $i$ \\
$n_{i,s}^{\text{total}}$ & Total number of sensors of type $s$ at site $i$ ($n_{i,s}^{\text{exist}} + n_{i,s}^{\text{add}}$) \\
$\beta_i$ & Binary variable indicating whether site $i$ is active (1) or inactive (0) \\
$\theta_{i,s,k}(t)$ & Minimum probability of detection for AAM aircraft $k$ by sensor $s$ at site $i$ at time $t$ \\
$R_u(t)$ & Server reliability at time $t$ \\
$\sigma$ & Detection threshold in robustness model \\
$\text{vert}_s$ & Number of vertical units required for sensor type $s$ \\
$\text{MaxVert}$ & Maximum vertical sensor units allowed per site \\
$\text{MaxSets}_s$ & Maximum allowed stack (copies) for sensor type $s$ at a site \\
$\text{cost}_s$ & Cost of deploying one unit of sensor type $s$ \\
$S$ & Set of all deployed sensors ($n_{i,s}^{\text{total}}$ for all $i,s$) \\
$\text{fail}[o,t]$ & Binary indicator if primary sensor $o$ is failed at time $t$ \\
$\text{fail\_times}[o]$ & Start time of failure for primary sensor $o$ \\
$\text{repair\_times}[o]$ & Repair duration for primary sensor $o$ \\
$\text{travel\_time}[h,o]$ & Time required for a backup sensor to travel from hub $h$ to primary sensor $o$ \\
$\text{dispatch}[b,h,o,t]$ & Binary decision: 1 if backup sensor $b$ is dispatched from hub $h$ to \\
& primary sensor $o$ at time $t$ \\
$\text{active}[b,o,t]$ & Binary decision: 1 if backup sensor $b$ is active at primary sensor $o$ \\
& at time $t$ \\
$\text{prob}_b$ & Detection probability of backup sensor $b$ \\
\end{longtable}

\section{Introduction}
Advanced Air Mobility (AAM) is a developing aviation concept that incorporates drones, including electric vertical takeoff and landing (eVTOL) aircraft and uncrewed aerial vehicles (UAVs), into low-altitude airspace to enable applications such as passenger transport, cargo delivery, and logistics services. AAM has the potential to provide substantial economic, social, and environmental benefits, particularly in Ohio \cite{Dulia2021AAM,DelRosario2021Infrastructure,Dulia2022OpenFramework,Calhoun2023OpenFramework}.
 The safe and efficient operation of AAM aircraft critically depends on continuous surveillance of low-altitude airspace, enabling functions such as low-altitude traffic management, detection and avoidance in Beyond Visual Line of Sight (BVLOS) operations, and counter-UAS (detect-and-intercept) capabilities.

A Surveillance for Advanced Air Mobility (SAM) system is envisioned as the cornerstone of this infrastructure. Such a system typically consists of a multi-type sensor network, data links, and a cloud-based information clearinghouse, together enabling detection, identification, and tracking of all AAM aircraft in the operating region. Although prior research has extensively explored sensor network design and nominal coverage optimization, most studies assume ideal operating conditions and do not account for real-world operational uncertainties. In practice, SAM systems are subject to sensor and communication link failures, adverse environmental conditions, high traffic density, and possible server outages, all of which can compromise system performance. Furthermore, in the event of a major failure or system-wide outage, SAM operations must be able to recover quickly through backup resources to ensure continuity of surveillance and airspace safety. These gaps highlight the need for a quantitative framework that integrates reliability, robustness, and resilience (3Rs) into the planning of the SAM system.

The concepts of 3Rs form a foundational framework for analyzing and improving the performance of complex engineered systems under uncertainty and disruptions. Reliability refers to the probability that a system performs its intended function successfully over a specified period of time under expected operating conditions \cite{ReliabilityTextbook,BarlowProschan,ModarresReliabilityEngineering}. Robustness describes the ability of a system to maintain acceptable performance despite uncertainties, variations, or partial degradations in system components or inputs \cite{BenTalRobustOptimization,HaleRobustSystems,Holland2003RobustnessComplexSystems}. In contrast, resilience captures the system’s ability to rapidly recover and restore functionality following significant disruptions or failures \cite{Bruneau2003Resilience,Francis2018ResilienceReview}. These three properties provide a unified framework for designing systems that perform reliably under normal conditions, remain stable under uncertainty, and recover effectively after failures.



In our previous work \cite{Dulia2024SAM}, we developed an optimization model that determines the optimal number and placement of sensors to design a surveillance network for monitoring AAM aircraft. However, that study did not explicitly incorporate the 3R principles of reliability, robustness, and resilience. A primary limitation preventing 3R-compatible design for AAM surveillance was that the formulation focused on optimal sensor placement within a simplified two-dimensional framework to achieve nominal geometric coverage. However, a 2D representation cannot adequately capture line-of-sight (LOS) constraints, which are critical in AAM surveillance. In realistic urban environments, sensor detection is strongly affected by elevation differences and obstructions caused by tall buildings, terrain, and other natural or human-made barriers, which may partially degrade or completely block visibility between sensors and aircraft. In addition, the previous formulation did not account for system component reliability, including potential failures of sensors, communication links, and servers, all of which directly influence effective detection and tracking performance beyond geometric coverage. Moreover, sensor tracking capability is inherently dependent on aircraft density and traffic load. In practical AAM operations, traffic is spatially and temporally heterogeneous, with certain regions experiencing high-density operations while others remain sparse. This variability directly affects sensor workload and maximum tracking capacity, which was not captured in the previous study due to the assumption of a uniform spatial distribution of aircraft. As a result, the simplified formulation does not adequately represent operational stress conditions such as sudden traffic surges or localized congestion, which are essential for evaluating system robustness and resilience.

These limitations indicate that the previous study is not fully compatible with the 3R framework and therefore requires further refinement. To address this gap, this study incorporates 3D spatial modeling with LOS constraints, modeling of component reliability, and realistic AAM flight schedules that capture variable traffic demand. Together, these enhancements enable a more comprehensive and realistic evaluation of surveillance performance under normal, degraded, and high-demand conditions, thereby supporting 3R-oriented design of AAM surveillance systems. As part of the 3R framework, the reliability model determines the optimal placement and selection of sensors to meet minimum detection thresholds under nominal operating conditions. The robustness model, currently under development, extends the framework to maintain acceptable system performance under off-nominal conditions such as partial sensor failures, adverse weather, or high-density traffic by identifying additional sensor placements and optimal data transmission strategies. Finally, the resilience model, also under development, plans backup sensor allocation and dispatch strategies to rapidly restore coverage following major failures, thereby minimizing system downtime and ensuring the safe resumption of AAM operations.

\section{Methodology}

\begin{figure}[th] 
    \centering
    \includegraphics[width=17cm,height=10cm]{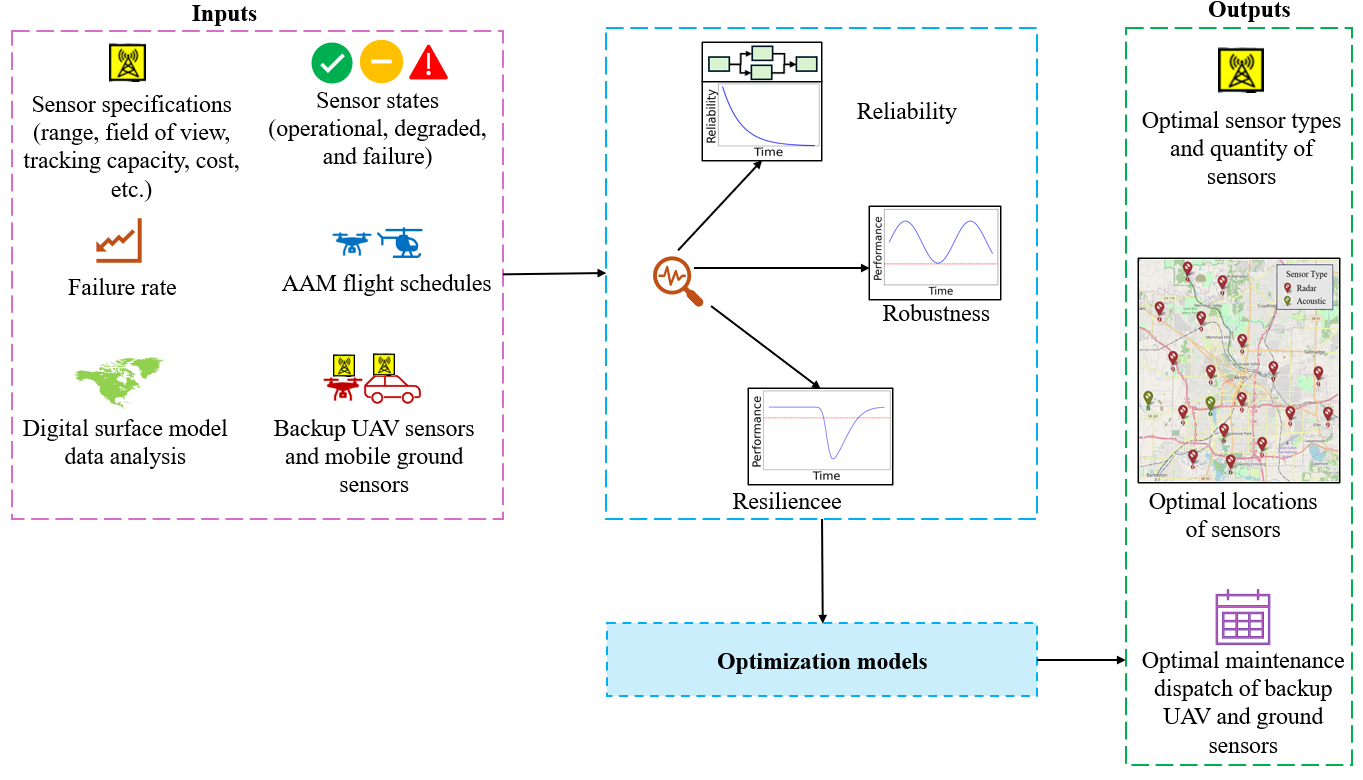}
    \caption{Overview of our approach for developing the 3R framework for the SAM system.} 
    \label{fig_i_80_output2}
\end{figure}

In this study, we follow a structured approach to develop the 3R models for the SAM system. Figure~\ref{fig_i_80_output2} illustrates an overview of the proposed 3R modeling framework for the SAM system. The framework integrates diverse inputs, including sensor specifications, SAM component states and their associated failure rates, projected AAM flight schedules, spatial and terrain data, and other relevant operational parameters. These inputs are then used in optimization models developed around the three core components of the 3R framework: reliability, robustness, and resilience. The outputs include the optimal sensor types and quantities, their placement along AAM corridors, and operational strategies to maintain continuous surveillance coverage under both nominal and off-nominal conditions.

The surveillance environment consists of a set of candidate sites $\mathcal{I}$, each capable of hosting multiple sensor types from the set $\mathcal{S}$. A set of AAM aircraft $\mathcal{K}$ operates within the surveillance region and must be detected over the analysis time horizon $\mathcal{T}$. The mathematical formulations of the three 3R components are presented next.

\subsection{Reliability Model} 

The reliability model represents the first component of the 3R framework and determines the optimal placement of heterogeneous sensors such that every AAM aircraft is reliably detected and its data are successfully transmitted to the data server. The main decision variables are $n_{i}^{s}$, representing the number of sensor sets of type $s$ installed at site $i$, and $\beta_i$, a binary indicator that equals 1 if site $i$ is activated and 0 otherwise. The mathematical formulation of the reliability model is presented in the following equations.

The reliability model begins by minimizing the total installation cost of deploying heterogeneous sensors across candidate locations, as formulated in Equation \ref{eq:obj}. The objective function $Z$ computes the overall deployment cost by combining the installation cost of each sensor type $c_s$, the number of sensors contained in one set $C_s$, and the number of installed sensor sets $n_i^s$ at site  $i$. Since the operational reliability of the surveillance system changes over time, Equation \ref{eq:reliability} models the time-dependent reliability of the key system components—sensors, communication links, and the central server—using exponential decay functions defined by their failure rates $\lambda_s$, $\lambda_\ell$, and $\lambda_u$, which determine the reliability values $\rho_s(t)$, $\rho_\ell(t)$, and $\rho_u(t)$ at time $t$. To determine whether an AAM aircraft can be observed by a sensor, Equation \ref{eq:range} introduces a binary range indicator $\chi_{i,s,k}(t)$ that checks whether the three-dimensional distance $d_{i,k}(t)$ between candidate site $i$ and aircraft $k$ at time $t$ lies within the effective sensing range $r_s$ of sensor type $s$; this feasibility check is implemented through the \textit{RangeCheck} procedure (Algorithm \ref{alg:rangecheck}). Next, Equation \ref{eq:intrinsic} defines the intrinsic detection probability $p_{i,s,k}(t)$ by combining sensor reliability $\rho_s(t)$, range feasibility $\chi_{i,s,k}(t)$, and the LOS quality factor $\ell_{i,s,k}(t)$, which represents the fraction of unobstructed visibility between the sensor and the aircraft. 

\begin{equation}
\min Z = \sum_{i \in \mathcal{I}} \sum_{s \in \mathcal{S}} c_{s} C_{s} n_{i}^{s}
\label{eq:obj}
\end{equation}

\begin{equation}
\rho_{s}(t) = e^{-\lambda_{s} t}, \quad
\rho_{\ell}(t) = e^{-\lambda_{\ell} t}, \quad
\rho_{u}(t) = e^{-\lambda_{u} t}, \quad \forall t \in \mathcal{T}
\label{eq:reliability}
\end{equation}

\begin{equation}
\chi_{i,s,k}(t) =
\begin{cases}
1, & d_{i,k}(t) \le r_{s}, \quad \forall i \in \mathcal{I}, \forall s \in \mathcal{S}, \forall k \in \mathcal{K}, \forall t \in \mathcal{T}\\[6pt]
0, & \text{otherwise}
\end{cases}
\label{eq:range}
\end{equation}

\begin{algorithm}[h]
\caption{Range Check Function}
\label{alg:rangecheck}
\begin{algorithmic}[1] 
\Function{RangeCheck}{S, D, r, e, a, b}
    \State \textbf{Input:} 
        \State \quad S: 3D position of sensor
        \State \quad D: 3D position of AAM aircraft
        \State \quad r: detection radius of sensor
        \State \quad e: radius tolerance for boundary decay
        \State \quad a, b: decay parameters
    \State Compute Euclidean distance: $d \gets \|S - D\|$
    \If{$d \le r - e$} 
        \State \Return 1.0 \Comment{AAM aircraft well inside detection radius}
    \ElsIf{$r - e < d < r + e$} 
        \State \Return $\exp\Big(-a \,(d - (r - e))^b\Big)$ \Comment{Smooth decay near boundary}
    \Else
        \State \Return 0.0 \Comment{AAM aircraft outside detection radius}
    \EndIf
\EndFunction
\end{algorithmic}
\end{algorithm}

\begin{equation}
p_{i,s,k}(t) = \rho_{s}(t)\, \chi_{i,s,k}(t)\, \ell_{i,s,k}(t), \quad \forall i \in \mathcal{I}, \forall s \in \mathcal{S}, \forall k \in \mathcal{K}, \forall t \in \mathcal{T}
\label{eq:intrinsic}
\end{equation}

The LOS constraints are considered to capture visibility conditions between sensors and aircraft. The LOS factor $\ell_{i,s,k}(t)$ is computed using a geometric LOS formulation in Algorithm \ref{alg:geolos} for electromagnetic wave–based sensing systems, including radar, RF sensors, ADS-B and Remote ID receivers, and optical cameras. The specific sensor types used in the model are further detailed in the sensor specification section. In this context, electromagnetic waves refer to radio-frequency signals used by radar, RF sensors, ADS-B, and Remote ID receivers, as well as light waves used by optical cameras, all of which can be obstructed by terrain or physical structures such as buildings. The model determines whether there is a clear, unobstructed path between the sensor and the aircraft by checking whether such obstacles intersect the straight line connecting them. If any obstruction exists along this path, the LOS is reduced accordingly. For radar, RF sensors, ADS-B, and Remote ID receivers, this LOS represents the availability of a clear signal propagation path. For optical cameras, it represents direct visual visibility. Although the physical sensing mechanisms differ, all sensor types are modeled using the same geometric visibility principle in this study. This LOS formulation is based on classical viewshed analysis and terrain-based visibility computation methods used in geographic information systems and computational geometry. In particular, the concept of determining visibility by testing whether terrain intersects the straight line between an observer and a target follows standard viewshed algorithms widely used in GIS-based spatial analysis \cite{reddy2018gis,gomarasca2009basic}. The projection of terrain points onto the line segment and evaluation of obstruction conditions is consistent with ray-casting and line-segment intersection techniques commonly used in computational geometry and terrain modeling \cite{deberg2008computational}. These methods are also widely applied in DSM-based visibility analysis, where obstruction is determined by comparing terrain elevation profiles against interpolated LOS paths \cite{duckham2023gis}.

\begin{algorithm}[h]
\caption{Geometric LOS Function}
\label{alg:geolos}
\begin{algorithmic}[1]
\Function{LOS}{S, D, T, R, Z}
    \State \textbf{Input:} 
        \State \quad S: 3D position of sensor
        \State \quad D: 3D position of AAM aircraft
        \State \quad T: terrain points
        \State \quad R: horizontal cylinder radius tolerance
        \State \quad Z: vertical tolerance
    \If{$S = D$} 
        \State \Return 1.0
    \EndIf
    \State Project positions to XY-plane: $S_{xy}, D_{xy}$
    \State Compute vector $\vec{V} = D_{xy} - S_{xy}$ and length $L = \|\vec{V}\|$
    \State Project terrain points onto line: $\text{proj} = (T_{xy}-S_{xy}) \cdot \vec{V}/L$
    \State Compute horizontal distance: $d_\text{line} = \|\text{T}_{xy} - \text{closest\_point\_on\_line}\|$
    \State Mask points satisfying $0 \le \text{proj} \le L$ and $d_\text{line} \le R$
    \If{no points to check} \State \Return 1.0 \EndIf
    \State Compute fraction along line: $\tau = \text{proj}/L$
    \State Interpolate line heights: $h = S_z + \tau (D_z - S_z)$
    \State Determine blocked points: $T_z \ge h - Z$
    \State Compute geometric LOS: $\ell_\text{geometric} = 1 - \frac{\# \text{blocked points}}{\# \text{points to check}}$
    \State \Return $\ell_\text{geometric}$
\EndFunction
\end{algorithmic}
\end{algorithm}

While the above LOS formulation applies to electromagnetic sensing modalities, acoustic sensors require a distinct propagation model due to their fundamentally different wave physics. The LOS factor $\ell_{i,s,k}(t)$ is computed using Algorithm \ref{alg:acoustillos} for acoustic sensors. Acoustic sensors operate based on sound wave propagation, where aircraft are detected through pressure waves traveling through air. The acoustic LOS algorithm is constructed by integrating well-established principles from classical acoustics, outdoor sound propagation theory, and terrain-dependent attenuation models. The spherical wave propagation and distance-dependent amplitude decay are grounded in fundamental acoustic theory, where sound intensity decreases with distance and is described using standard wavenumber formulations \cite{pierce2019acoustics,kinsler2000fundamentals}. The geometric spreading and reference attenuation components follow the ISO 9613-2 standard, which is widely used for modeling outdoor sound propagation and environmental noise attenuation \cite{iso9613}. Terrain-induced effects, including diffraction, reflection, and shielding caused by natural and built environments, are incorporated using established computational acoustics and outdoor propagation frameworks \cite{salomons2001computational,attenborough2007predicting,gkanos2025numerical}.

\begin{algorithm}[h]
\caption{Acoustic LOS Function}
\label{alg:acoustillos}
\begin{algorithmic}[1]
\Function{LOS\_acoustic}{S, D, T, C, f, c, R, Z, A}
    \State \textbf{Input:} 
        \State \quad S: 3D position of sensor
        \State \quad D: 3D position of AAM aircraft
        \State \quad T: terrain points 
        \State \quad C: terrain types for acoustic attenuation
        \State \quad f: acoustic frequency
        \State \quad c: sound speed
        \State \quad R: horizontal cylinder radius
        \State \quad Z: vertical tolerance
        \State \quad A: acoustic amplitude coefficient
    \State Compute direct distance: $r = \|S - D\|$
    \If{$r < \delta$} \State \Return 1.0 \EndIf
    \State Compute wavenumber $k = 2 \pi f / c$
    \State Compute direct wave amplitude $\phi = A / r \cdot \exp(-i k r)$
    \State Identify nearby terrain points within $R$ of $D$ for reflection
    \State Compute geometric attenuation $G$ using LOS
    \State Compute barrier attenuation $B$ using terrain classes
    \State Total attenuation: $A_\text{tot} = G + B$
    \State Compute acoustic LOS: $\ell_\text{acoustic} = \frac{1}{1 + 10^{A_\text{tot}/20}}$
    \State \Return $\ell_\text{acoustic}$
\EndFunction
\end{algorithmic}
\end{algorithm}

Once detection probabilities are established for all sensors, the model extends to account for communication reliability and successful data transmission. Once an aircraft is detectable, Equation \ref{eq:transmission} extends the detection process to include communication reliability by defining the end-to-end probability $q_{i,s,k}(t)$ that a detected signal is successfully transmitted to the data server, which depends on both the intrinsic detection probability $p_{i,s,k}(t)$ and the communication link reliability $\rho_\ell(t)$. Equation \ref{eq:missed} then defines the missed detection probability $m_{i,s,k}(t)$ as the complement of the successful transmission probability, while enforcing a lower bound $\varepsilon$ to maintain numerical stability in the optimization. To track which aircraft are effectively detected and reported, Equation \ref{eq:coverage} defines the set $\alpha_{t,i,s}$ consisting of AAM aircraft $k$ that achieve a positive end-to-end detection probability $q_{i,s,k}(t)$ from sensor type $s$ at site $i$ during time step $t$. Using this information, Equation \ref{eq:reliability-constraint} enforces the overall system reliability requirement by aggregating the missed detection probabilities from all deployed sensors $n_i^s$ and ensuring that the combined detection performance satisfies the minimum required reliability threshold $H$, while also accounting for the reliability of the data server $\rho_u(t)$. Finally, practical deployment constraints are incorporated to ensure feasible infrastructure planning: Equation \ref{eq:capacity} restricts the total number of installed sensors at each candidate site according to the site capacity $M$ and the site activation decision variable $\beta_i$, Equation \ref{eq:linking} links sensor deployment to site activation while limiting the number of sensor sets of each type using $M_s$, and Equation \ref{eq:min-deploy} ensures that any activated site must host at least one sensor set, preventing the selection of empty sites.

\begin{equation}
q_{i,s,k}(t) = p_{i,s,k}(t)\, \rho_{\ell}(t), \quad \forall i \in \mathcal{I}, \forall s \in \mathcal{S}, \forall k \in \mathcal{K}, \forall t \in \mathcal{T}
\label{eq:transmission}
\end{equation}

\begin{equation}
m_{i,s,k}(t) = 1 - q_{i,s,k}(t), 
\qquad m_{i,s,k}(t) \ge \varepsilon, \quad \forall i \in \mathcal{I}, \forall s \in \mathcal{S}, \forall k \in \mathcal{K}, \forall t \in \mathcal{T}
\label{eq:missed}
\end{equation}

\begin{equation}
\alpha_{t,i,s} = \{ k \in \mathcal{K} : q_{i,s,k}(t) > 0 \}, \quad \forall i \in \mathcal{I}, \forall s \in \mathcal{S}, \forall k \in \mathcal{K}, \forall t \in \mathcal{T}
\label{eq:coverage}
\end{equation}

\begin{equation}
\sum_{\substack{i \in \mathcal{I},\, s \in \mathcal{S} \\ k \in \alpha_i^s(t)}} 
n_i^s \ln\!\big(m_{i,s,k}(t)\big)
\le 
\ln\!\Big(1 - \frac{H}{\rho_u(t)}\Big), \quad \forall t \in \mathcal{T}
\label{eq:reliability-constraint}
\end{equation}

\begin{equation}
\sum_{s \in \mathcal{S}} C_{s} n_{i}^{s} \le M \beta_{i}, \quad \forall i \in \mathcal{I}
\label{eq:capacity}
\end{equation}

\begin{equation}
n_{i}^{s} \le M_{s} \beta_{i},
\qquad \forall i \in \mathcal{I},\, \forall s \in \mathcal{S}, 
\label{eq:linking}
\end{equation}

\begin{equation}
\sum_{s \in \mathcal{S}} n_{i}^{s} \ge \beta_{i}, \quad \forall i \in \mathcal{I}
\label{eq:min-deploy}
\end{equation}

\newpage 

\subsection{Robustness Model}

The second model of the 3R models is the robustness model, which examines how the sensor network adjusts to sudden increases in AAM aircraft traffic density caused by rapidly changing weather, emergency AAM operations, or congestion in certain airspace corridors. Under these circumstances, the original network design—while optimal under normal conditions—may not provide sufficient coverage or capacity to maintain reliable detection of AAM aircraft. Thus, the purpose of this model is to determine the number and placement of additional sensors necessary to maintain a minimum detection threshold. The core idea is to retain all existing optimal sensors determined by the reliability model and augment them only when required. The combined sensor network, consisting of existing and new sensors, must deliver consistent performance even when the number of AAM aircraft rises dramatically in certain regions and time periods.

The following equations define the constraints and objective of the robustness model optimization framework used to determine the deployment of additional sensors while maintaining the required detection performance under high AAM traffic density. Equation \ref{robA_1} defines the total number of sensors $n_{i,s}^{\text{total}}$ of type $s$ installed at site $i$ as the sum of the existing sensors $n_{i,s}^{\text{exist}}$ obtained from the reliability model and the additional sensors $n_{i,s}^{\text{add}}$ determined by the optimization model. Based on these total sensors, Equation \ref{robA_2} enforces the coverage requirement by ensuring that the combined detection capability of all deployed sensors satisfies the minimum detection threshold $\sigma$ for each AAM aircraft $k$ at time $t$, where the detection contribution of each sensor is represented through the minimum probability of detection $\theta_{i,s,k}(t)$ and adjusted by the server reliability $R_u(t)$. To ensure feasible infrastructure deployment, Equation \ref{robA_3} limits the total vertical stacking of sensors at each site by constraining the sum of vertical units $\text{vert}_s$ required by sensors of type $s$ multiplied by their total counts $n_{i,s}^{\text{total}}$, ensuring that it does not exceed the maximum vertical capacity $\text{MaxVert}$ when the site activation variable $\beta_i$ indicates that site $i$ is active. Similarly, Equation \ref{robA_4} restricts the number of sensors of each type installed at a site by ensuring that the total number $n_{i,s}^{\text{total}}$ does not exceed the maximum allowable stack limit $\text{MaxSets}_s$ for sensor type $s$, while linking this constraint to the site activation decision $\beta_i$. Equation \ref{robA_5} guarantees that every activated site contains at least one sensor by requiring that the sum of all installed sensors across sensor types $\sum_{s \in \mathcal{S}} n_{i,s}^{\text{total}}$ is at least one whenever $\beta_i=1$, thereby preventing the selection of empty sites. The binary nature of the site activation decision is defined in Equation \ref{robA_6}, where $\beta_i$ takes a value of 1 if candidate site $i$ is activated and 0 otherwise. Equation \ref{robA_7} specifies that the number of additional sensors $n_{i,s}^{\text{add}}$ is a non-negative integer decision variable for each site $i$ and sensor type $s$, ensuring that sensor deployment decisions remain practically feasible. Equation \ref{robA_8} treats the number of existing sensors $n_{i,s}^{\text{exist}}$ obtained from the reliability model as fixed parameters, meaning that the robustness optimization only determines the number of additional sensors to deploy. Equation \ref{robA_9} then defines the set containing all deployed sensors across sites and sensor types, represented by the total sensor counts $n_{i,s}^{\text{total}}$, which is used for subsequent coverage and network calculations. Finally, the objective function in Equation \ref{robA_10} minimizes the total cost of deploying additional sensors by summing the deployment cost $\text{cost}_s$ of each sensor type multiplied by the required vertical units $\text{vert}_s$ and the number of additional sensors $n_{i,s}^{\text{add}}$ installed across all candidate sites $i$ and sensor types $s$. Together, these equations ensure that the sensor network maintains the minimum detection threshold $\sigma$ under high AAM aircraft density while respecting site activation decisions $\beta_i$, sensor stacking limits $\text{MaxVert}$ and $\text{MaxSets}_s$, and overall cost efficiency.

\begin{equation} \label{robA_1}
n_{i,s}^{\text{total}} = n_{i,s}^{\text{exist}} + n_{i,s}^{\text{add}}, \quad \forall i \in \mathcal{I}, \; s \in \mathcal{S}
\end{equation}

\begin{equation} \label{robA_2}
\sum_{i \in \mathcal{I}} \sum_{s \in \mathcal{S}} n_{i,s}^{\text{total}} \, \log(\theta_{i,s,k}(t)) 
\le \log\Big(1 - \frac{\sigma}{R_u(t)}\Big), 
\quad \forall k \in \mathcal{K}, \; t \in \mathcal{T}
\end{equation}

\begin{equation} \label{robA_3}
\sum_{s \in \mathcal{S}} \text{vert}_s \, n_{i,s}^{\text{total}} \le \text{MaxVert} \cdot \beta_i, 
\quad \forall i \in \mathcal{I}
\end{equation}

\begin{equation} \label{robA_4}
n_{i,s}^{\text{total}} \le \text{MaxSets}_s \cdot \beta_i, 
\quad \forall i \in \mathcal{I}, \; s \in \mathcal{S}
\end{equation}

\begin{equation} \label{robA_5}
\sum_{s \in \mathcal{S}} n_{i,s}^{\text{total}} \ge \beta_i, 
\quad \forall i \in \mathcal{I}
\end{equation}

\begin{equation} \label{robA_6}
\beta_i \in \{0,1\}, 
\quad \forall i \in \mathcal{I}
\end{equation}

\begin{equation} \label{robA_7}
n_{i,s}^{\text{add}} \ge 0, \quad n_{i,s}^{\text{add}} \in \mathbb{Z}, 
\quad \forall i \in \mathcal{I}, \; s \in \mathcal{S}
\end{equation}

\begin{equation} \label{robA_8}
n_{i,s}^{\text{exist}} = \text{given from reliability model}, 
\quad \forall i \in \mathcal{I}, \; s \in \mathcal{S}
\end{equation}

\begin{equation} \label{robA_9}
S = \{ n_{i,s}^{\text{total}} \,|\, i \in \mathcal{I}, s \in \mathcal{S} \}
\end{equation}

\begin{equation} \label{robA_10}
\min \sum_{i \in \mathcal{I}} \sum_{s \in \mathcal{S}} \text{cost}_s \cdot \text{vert}_s \cdot n_{i,s}^{\text{add}}
\end{equation}

\subsection{Resiliency Model}

While the robustness model looks at whether the system can continue working under off-nominal conditions and still meets a minimum acceptable performance level, the resilience model looks at how the system restores when performance drops below that level because of sensor failures. The resiliency model provides solutions to help the surveillance system keep operating even when one or more primary sensors fail. Its main purpose is to maintain continuous monitoring while failed sensors are repaired, minimizing gaps in coverage, and maintaining overall system performance. From our study in Task 1, resilience in a surveillance system is defined as its ability to quickly recover from failures or disruptions. The primary metric used to evaluate resilience is the recovery time, which measures how long it takes to restore a primary sensor from a failed state back to full operational status.

To achieve resiliency, the model considers backup sensors that temporarily assume the monitoring responsibilities of failed primary sensors. Consider a scenario in which primary sensors monitor AAM aircraft within a specific region. The number and placement of these primary sensors is determined using reliability and robustness models, which identify optimal sensor locations and coverage requirements. When a primary sensor fails, its performance decreases, and recovery procedures are initiated. During this recovery time, backup sensors are dispatched to provide temporary surveillance, ensuring continuous coverage of the SAM system. There are two types of backup sensors: i) UAV-mounted sensors, which can fly to the location of the failed sensor, and ii) ground vehicle–mounted sensors, which reach the location by road. Backup sensors are stored at pre-designated hubs and are deployed only when a primary sensor fails. Once the failed primary sensor is restored to normal operation, the backup sensor returns to its hub, ready for future deployments. This model allows the surveillance system to maintain operational capability during failures, reduce coverage gaps, and enhance overall resilience through proactive backup sensor management.

The following equations define the constraints and objectives of the resiliency model, which determines how backup sensors are dispatched and activated to maintain detection capability when primary sensors fail. Equation \ref{res_R_1} defines the activation feasibility of backup sensors by ensuring that a backup sensor $b \in \mathcal{B}$ can be active at a primary sensor location $o \in \mathcal{O}$ at time $t \in \mathcal{T}$ only if it has previously been dispatched from a candidate hub $h \in \mathcal{H}$, as represented by the binary dispatch decision variable $\text{dispatch}_{b,h,o,t}$. Equation \ref{res_R_2} further enforces this logic by setting the activation variable $\text{active}_{b,o,t}$ to zero when no valid dispatch window exists for backup sensor $b$ to reach primary sensor $o$ at time $t$, such as when the primary sensor has not failed or when the travel time $\text{travel\_time}[h,o]$ makes the dispatch infeasible. To prevent operational conflicts, Equation \ref{res_R_3} ensures that a backup sensor $b$ cannot be dispatched to another failed primary sensor $o_2$ during the repair period of a previously assigned failure at primary sensor $o_1$, where repair durations are determined by the failure start time $\text{fail\_times}[o]$ and repair duration $\text{repair\_times}[o]$. Equation \ref{res_R_4} further guarantees operational feasibility by restricting each backup sensor $b$ to be active at at most one primary sensor location at any time $t$, ensuring that a single backup platform cannot simultaneously cover multiple failures. Equation \ref{res_R_5} limits redundant assignments by ensuring that each backup sensor $b$ can be dispatched at most once to address a specific primary sensor failure at location $o$, within the failure period defined by $\text{fail}[o,t]$. To guarantee service continuity, Equation \ref{res_R_6} requires that at least one backup sensor from the set $B$ must be dispatched from some hub $h$ to each failed primary sensor $o$ during its failure period, thereby ensuring minimum coverage restoration. Finally, Equation \ref{res_R_7} defines the objective function of the resiliency model, which minimizes the total travel time required to dispatch backup sensors from hubs to failed primary sensors using the parameter $\text{travel\_time}[h,o]$, while simultaneously maximizing their effective detection contribution by rewarding active backup sensors through their detection probability $\text{prob}_b$. Together, these equations ensure that backup sensors are dispatched efficiently from candidate hubs, avoid operational conflicts, and maintain detection coverage during primary sensor failures.

\begin{equation}\label{res_R_1}
\text{active}_{b,o,t} \leq \sum_{h \in H} \sum_{t_{dispatch} \le t} \text{dispatch}_{b,h,o,t_{dispatch}}, \quad \forall b \in \mathcal{B}, o \in \mathcal{O}, t \in \mathcal{T}
\end{equation}

\begin{equation}\label{res_R_2}
\text{active}_{b,o,t} = 0, \quad \text{if no valid dispatch window exists for } b, o, t
\end{equation}

\begin{equation}\label{res_R_3}
\sum_{h \in H} \text{dispatch}_{b,h,o_2,t_2} \le 1 - \sum_{h \in \mathcal{H}} \sum_{t_w} \text{dispatch}_{b,h,o_1,t_w}, \quad \forall b \in \mathcal{B}, o_1 \neq o_2, t_2 \in \text{repair\_period}(o_1)
\end{equation}

\begin{equation}\label{res_R_4}
\sum_{o \in \mathcal{O}} \text{active}_{b,o,t} \le 1, \quad \forall b \in \mathcal{B}, t \in \mathcal{T}
\end{equation}

\begin{equation}\label{res_R_5}
\sum_{h \in \mathcal{H}} \sum_{t \in \text{fail\_period}(o)} \text{dispatch}_{b,h,o,t} \le 1, \quad \forall b \in \mathcal{B}, o \in \mathcal{O}
\end{equation}

\begin{equation}\label{res_R_6}
\sum_{b \in \mathcal{B}} \sum_{h \in \mathcal{H}} \sum_{t \in \text{fail\_period}(o)} \text{dispatch}_{b,h,o,t} \ge 1, \quad \forall o \in \mathcal{O}
\end{equation}

\begin{equation}\label{res_R_7}
\min \sum_{b \in \mathcal{B}} \sum_{h \in \mathcal{H}} \sum_{o \in \mathcal{O}} \sum_{t \in \mathcal{T}} \text{travel\_time}_{h,o} \cdot \text{dispatch}_{b,h,o,t} - \sum_{b \in B} \sum_{o \in O} \sum_{t \in T} \text{prob}_b \cdot \text{active}_{b,o,t}
\end{equation}

\section{Experimental Setup}

This section describes the datasets, modeling assumptions, and parameters used to set the experimental environment to evaluate the 3R models. The experimental setup incorporates geographic terrain data, projected flight demand, sensor system characteristics, communication infrastructure parameters, and other modeling inputs.

\subsection{Digital Surface Model Data Analysis}

To demonstrate the performance of the reliability model, we applied the model to the five major corridors (I-70, I-71, I-75, I-80, and US-33) identified for the development of the SAM sensor network, as mentioned in the data collection section. We used DSM data collected for the geographic regions surrounding each corridor. DSM data represent high-resolution elevation information that includes the terrain surface as well as above-ground features such as buildings, vegetation, and other man-made structures. The DSM data were obtained from publicly available LiDAR datasets provided by the United States Geological Survey (USGS) through the National Map Downloader \cite{usgs_tnm}. These datasets provide high-resolution spatial information, including geographic coordinates, elevation values, and surface characteristics derived from thousands of LiDAR points across the study area.

\begin{figure}[tb!]
\centering
\includegraphics[width=8cm,height=8cm]{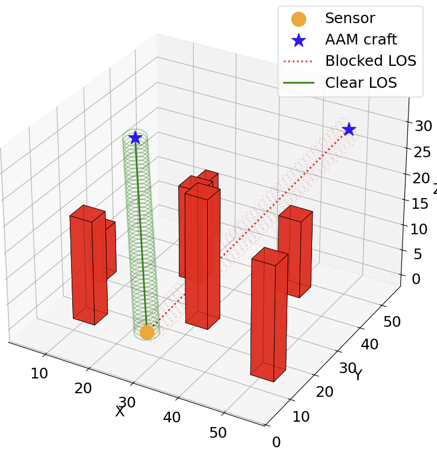}
\caption{Demonstration of LOS between a sensor and an AAM aircraft. Terrain features can block or attenuate sensor signals, affecting detection performance.}
\label{los}
\end{figure}

The inclusion of DSM data enables los analysis between sensors and AAM aircraft, which is critical for accurately modeling sensor detection performance. Terrain features, including hills, buildings, and vegetation, can obstruct or attenuate sensor signals, thereby reducing detection range and reliability. Figure~\ref{los} illustrates a representative LOS scenario between a sensor and an AAM aircraft, highlighting how terrain features can block or partially obscure the sensor’s field of view. Incorporating LOS effects ensures that the reliability model accounts for real-world constraints imposed by terrain and man-made structures, leading to more accurate assessments of surveillance coverage.

\begin{figure}[tb!]
\centering
\includegraphics[width=15cm,height=9cm]{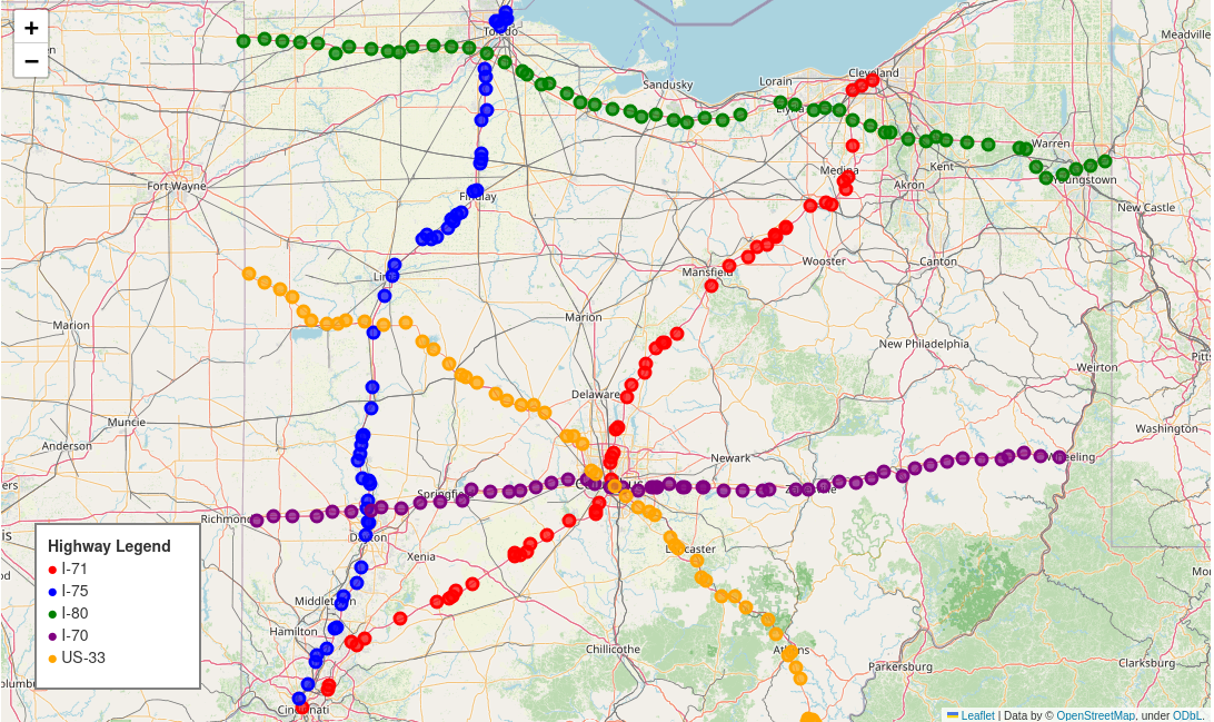}
\caption{Overview of the five Ohio corridors considered in the analysis.}
\label{fig_corridors}
\end{figure}

The downloaded LiDAR files were provided in compressed \texttt{.laz} format, which contains point-cloud information describing surface characteristics. These files were converted to CSV format to facilitate data processing. Each point includes longitude, latitude, altitude, and a terrain classification label based on the standard LAS classification scheme. The LAS classification system assigns numerical codes ranging from 1 to 18, where each value corresponds to a specific terrain type. However, the downloaded datasets for the selected corridors contained only three classification values after conversion: 1) Unclassified, 2) Ground (bare earth), and 5: High vegetation (trees and tall shrubs). Table~\ref{tab:dsm_sample} presents a representative sample of the DSM dataset for the I-70 corridor.

\begin{table}[h!]
\centering
\caption{Sample DSM Data from the I-70 Corridor}
\label{tab:dsm_sample}
\begin{tabular}{cccc}
\toprule
Longitude & Latitude & Altitude (ft) & Classification \\
\midrule
-81.02885608 & 40.06759534 & 1193.69 & 1 \\
-81.02883988 & 40.06836594 & 1191.92 & 2 \\
-81.02876972 & 40.06903201 & 1215.71 & 2 \\
-81.02855269 & 40.06974697 & 1215.40 & 2 \\
-81.02863442 & 40.07030448 & 1271.38 & 5 \\
-81.02838621 & 40.07062580 & 1278.62 & 5 \\
-81.02807263 & 40.07134216 & 1240.13 & 1 \\
-81.02805489 & 40.07181763 & 1256.13 & 2 \\
-81.02848141 & 40.07280399 & 1224.98 & 2 \\
\vdots & \vdots & \vdots & \vdots \\
\bottomrule
\end{tabular}
\end{table}

Across all datasets, a large proportion of LiDAR points were labeled as unclassified, indicating incomplete terrain labeling in the original survey data. To address this limitation, additional terrain classifications were inferred on the basis of elevation differences between LiDAR points by applying the following assumptions:

\begin{itemize}
\item Low vegetation (Class 3): height $<0.5$ m above ground level,
\item Medium vegetation (Class 4): vegetation between $0.5$ m and $5$ m, and
\item Buildings (Class 6): structures that exceed $5$ m in height.
\end{itemize}

\begin{figure}[th]
\centering
\includegraphics[width=15cm,height=6cm]{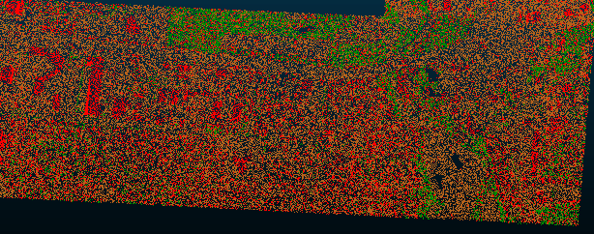}
\caption{Terrain visualization for a segment of the I-70 corridor derived from LiDAR data. Red represents buildings, green represents vegetation, and brown represents ground surfaces.}
\label{fig_i70cloud}
\end{figure}

Figure~\ref{fig_i70cloud} shows a small corridor segment of that corridor with different terrain classifications (red: buildings, green: vegetation, brown: ground). These variations in terrain were incorporated into the model to ensure that the results reflect the influence of terrain features on sensor placement.

\subsection{AAM Flight Schedule Generation}

the flight schedules were generated for three representative AAM use cases: regional passenger transportation (“Air Metro”), emergency medical transportation, and cargo delivery operations. The projected economic demand for these operations was obtained from an economic impact report published by the Ohio Department of Transportation \cite{rosario2021infrastructure}, which provides forecasts for advanced autonomous aircraft technologies along major Ohio transportation corridors through the year 2045. Based on these projections, the expected number of daily AAM flights for the year 2030 was estimated for each corridor. The resulting flight demand values are shown in Table~\ref{tab:daily_trips}.

\begin{table}[h!]
\centering
\caption{Estimated daily AAM trips per corridor}
\label{tab:daily_trips}
\begin{tabular}{lccc}
\toprule
Flight Corridor & Air Metro & Emergency & Cargo \\
\midrule
I-70 & 53 & 109 & 21 \\
I-71 & 162 & 272 & 42 \\
I-75 & 87 & 133 & 13 \\
I-80 & 69 & 146 & 26 \\
US-33 & 39 & 77 & 2 \\
\bottomrule
\end{tabular}
\end{table}

The daily flight counts were distributed temporally using realistic demand patterns. Emergency transportation flights were modeled using a bimodal distribution with peak activity occurring in the late morning and evening hours \cite{cantwell2015ems}. The demand for passenger transportation followed a trimodal distribution corresponding to the morning, noon, and evening commute periods. Cargo operations were modeled using a bimodal business-day distribution consistent with freight activity patterns reported in Ohio transportation studies \cite{odot2023freight,odot2018traffic}. Using these distributions, hourly flight schedules were generated between 9:00 AM and 6:00 PM for each corridor, and geographic trajectory information was integrated using the DSM data described above. Each aircraft trajectory consists of latitude, longitude, and altitude coordinates representing the aircraft position at each time step. Table~\ref{tab:flight_schedule} presents a sample flight schedule for the I-80 corridor. All AAM aircraft were assumed to operate at a constant cruising altitude of 400 feet.

\begin{table}[h!]
\centering
\caption{Sample AAM flight schedule for the I-80 corridor} \label{tab:flight_schedule}
\label{tab:aam_trajectories}
\begin{tabular}{p{1.5cm} p{3.5cm} p{3.5cm} p{3.5cm} p{1cm}}
\toprule
\textbf{Time} & \textbf{Air Metro Craft 1} & \textbf{Air Metro Craft 2} & \textbf{Air Metro Craft 3} & $\cdots$ \\
\midrule
9:00  & - & -84.0259590, 39.8637967, 400.00 & -81.2859993, 40.0597580, 400.00 & $\cdots$ \\

9:10  & - & -82.4682418, 39.9424398, 400.00 & -81.6445061, 40.0080039, 400.00 & $\cdots$ \\

9:20  & -81.2859993, 40.0597580, 400.00 & -81.2859993, 40.0597580, 400.00 & -81.9236234, 39.9618003, 400.00 & $\cdots$ \\

9:30  & -81.6471704, 40.0079980, 400.00 & - & -82.2264278, 39.9490844, 400.00 & $\cdots$ \\

9:40  & -81.9236234, 39.9618003, 400.00 & - & -82.5579169, 39.9425168, 400.00 & $\cdots$ \\

9:50  & -82.2264278, 39.9490844, 400.00 & -81.2859993, 40.0597580, 400.00 & -82.7977297, 39.9351986, 400.00 & $\cdots$ \\



$\vdots$ & $\vdots$ & $\vdots$ & $\vdots$ & $\ddots$ \\
\bottomrule
\end{tabular}
\end{table}

\subsection{Acoustic Propagation Modeling}

To accurately model acoustic sensor performance, sound propagation was considered across different types of terrain because the AAM corridors span diverse landscapes throughout Ohio. The characteristics of the terrain influence the reflection, absorption, and attenuation of sound, which directly affect the detection capabilities of acoustic sensors.

The acoustic reflection factor used in the analysis was calculated using the following model:

\begin{equation}
Q = R_p + (1-R_p)F(\omega)
\end{equation}

where $Q$ represents the reflection factor, $R_p$ denotes the ground reflection coefficient that describes the proportion of sound energy reflected by the terrain surface, and $F(\omega)$ represents the boundary loss factor associated with the ground impedance and frequency-dependent acoustic effects.

Terrain-specific values of $R_p$ were collected from studies of outdoor sound propagation and international acoustic standards \cite{iso9613,hannah2007wind}. These values vary depending on the composition of the terrain and the moisture conditions. Table~\ref{tab:terrain_rp} summarizes the ranges of reflection coefficients used in the model.

\begin{table}[h!]
\centering
\caption{Terrain Reflection Coefficient Values}
\label{tab:terrain_rp}
\begin{tabular}{lc}
\toprule
Terrain Classification & Ground Reflection Coefficient ($R_p$) \\
\midrule
Ground & 0.5--0.7 \\
Low Vegetation ($<$0.5 m) & 0.2--0.4 \\
Medium Vegetation (0.5--5 m) & 0.1--0.3 \\
High Vegetation (trees) & 0.05--0.2 \\
Buildings & 0.9--1.0 \\
\bottomrule
\end{tabular}
\end{table}

The boundary loss factor $F(\omega)$ was obtained from UAV acoustic detection studies \cite{dumitrescu2020acoustic}. Based on these studies, values in the range $0.975 \leq F(\omega) \leq 0.999$ were used depending on the acoustic impedance of the ground surface. This range is appropriate because AAM aircraft typically produce acoustic signatures within the frequency range of 200 to 10,000 Hz.

\subsection{Characteristics of Primary Sensors} \label{sensor_spec}

To evaluate the performance of the 3R models, parameters for several selected sensor types were collected from commercially available systems. These sensors include radar, optical cameras, Remote ID receivers, acoustic sensors, RF sensors, and ADS-B receivers. For each sensor type, key characteristics were gathered from manufacturer specifications and technical documentation. The collected parameters include detection range, cost, field of view, tracking capacity (the maximum number of UAVs the sensor can handle simultaneously), and failure rate \cite{dulia20263r}. Table~\ref{tab:sensor_specs} summarizes the specifications of the selected sensor systems considered in this study.

\begin{table}[h!]
\centering
\caption{Specifications of selected primary sensors}
\label{tab:sensor_specs}
\begin{tabular}{p{1.2cm} p{3.5cm} p{1.5cm} p{1.5cm} p{1.5cm} p{2cm} p{2cm}}
\toprule
Sensor Type & Manufacturer / Model & Tracking Capacity & Cost (per Unit) & Detection Range (km) & Failure Rate (failures/hr) & Field of View \\
\midrule
Radar & Echodyne – Echoflight & 20 craft & \$20,000 & Small UAV: 0.75; large UAV: 1  & $1.203\times10^{-5}$ & $120^\circ \times 80^\circ$ \\

Optical Camera & Edge Autonomy – Octopus ISR E180 & $>$500 craft & \$231,000 & 1.4–12.4 & $1.538\times10^{-5}$ & $360^\circ \times 30^\circ$ \\

Remote ID & Dronetag – RIDER & 20 craft & \$1,099 & 5 & $2.405\times10^{-5}$ & $360^\circ \times 180^\circ$ \\

Acoustic & DroneShield – FarAlert & 5 craft & $\sim\$20,000$ & 1 & $1.203\times10^{-5}$ & $360^\circ \times 90^\circ$ \\

RF & Aaronia AG – AARTOS X2 DDS & Unlimited & $\sim\$20,000$ & 5 & $6.014\times10^{-6}$ & $360^\circ \times 180^\circ$ \\

ADS-B & uAvionix – pingUSB & $\sim$100 craft (estimated) & \$275 & 160.9 & $6.014\times10^{-6}$ & $360^\circ \times 180^\circ$ \\
\bottomrule
\end{tabular}
\end{table}

\subsection{Characteristics of Backup Sensors}

To maintain continuous surveillance coverage, backup sensors are incorporated into the resiliency model. In the event of a primary sensor failure, secondary sensor units are deployed to temporarily replace the failed sensor while repair or replacement activities are performed. Backup sensors are stored at designated hub locations throughout the corridor network. When a primary sensor failure occurs, the nearest hub dispatches a backup unit to restore monitoring capability in the affected area. After the primary sensor becomes operational again, the backup unit returns to its hub and remains available for future deployments. Seven types of backup sensors are incorporated into the experimental setup: six are UAV-mounted and one is ground-based. Table~\ref{tab:backup_sensors} summarizes their key characteristics, including sensor weight, UAV or vehicle compatibility, detection range, probability of detection, and the maximum number of aircraft tracks that can be monitored simultaneously.

\begin{table}[tb!h]
\centering
\caption{Backup sensor characteristics}
\label{tab:backup_sensors}
\begin{tabular}{p{2cm} p{2cm} p{2cm} p{2.5cm} p{2.5cm} p{1.8cm}}
\toprule
Sensor & Weight (with Payload) & Compatible UAVs / Vehicle & Detection Range & Detection Probability & Tracking Capacity \\
\midrule
Echodyne Echoflight (Radar) & 3.9 kg & Freefly Astro & Small UAV: 0.75 km; Large UAV: 1 km &  0.95  & 100 \\
Dronetag RIDER (Remote ID) & 3.17 kg & Skydio X10, Freefly Astro & 5 km & 1.0 & 20 \\
uAvionix pingUSB (ADS-B) & 3.11 kg & Skydio X10, Freefly Astro & 160.9 km & 1.0 & 100 \\
CRYSOUND CRY2626G (Acoustic) & 4.5 kg & Freefly Astro & 35 m & 0.6--0.85 & 3 \\
Aaronia ISOLOG 3D (RF) & 3.45 kg & Freefly Astro & 1--3 km (urban: 0.5--1 km) & 0.7--0.9 & 10 \\
Shield AI ViDAR (Optical) & 5.1 kg & Freefly Astro & 3--5 km & 0.9--0.96 & 20 \\
DroneShield DroneSentry-X Mk2 (RF/Acoustic) & 46 kg & Ground Vehicle & 15 km & 0.95 (open), 0.7 (hills/vegetation) & Unlimited \\
\bottomrule
\end{tabular}
\end{table}

\newpage

\section{Results}

This section presents the results derived from the 3R models: the reliability model, robustness models A and B, and the resiliency model. The reliability model determines the optimal sensor type, the number of sensors, and their placement under normal operating conditions to meet SAM’s minimum detection requirements. Robustness Model A determines the additional number and placement of sensors needed when AAM traffic increases because of external perturbations, while robustness model B identifies the optimal data transmission paths and also measures performance loss when certain sensors fail to transmit data because of the perturbations. When the performance loss measured by robustness model B exceeds an acceptable threshold, the resiliency model is designed to assess the system’s ability to restore SAM performance. The results from these models are presented in the following subsections.

\subsection{Reliability Model}

To determine the number and placement of sensors needed along each corridor under normal AAM traffic conditions, the reliability model requires three inputs: (1) the corridor of interest, (2) the minimum reliability threshold set by SAM, and (3) the subset of sensor types considered in the analysis. Using these inputs, the model identifies the optimal configuration of sensors that satisfies the reliability target at the minimum possible cost. We evaluated several test cases to understand how sensor type selection and reliability requirements influence the overall network design.

\begin{figure}[tb!h]
\centering
\includegraphics[width=10cm,height=7cm]{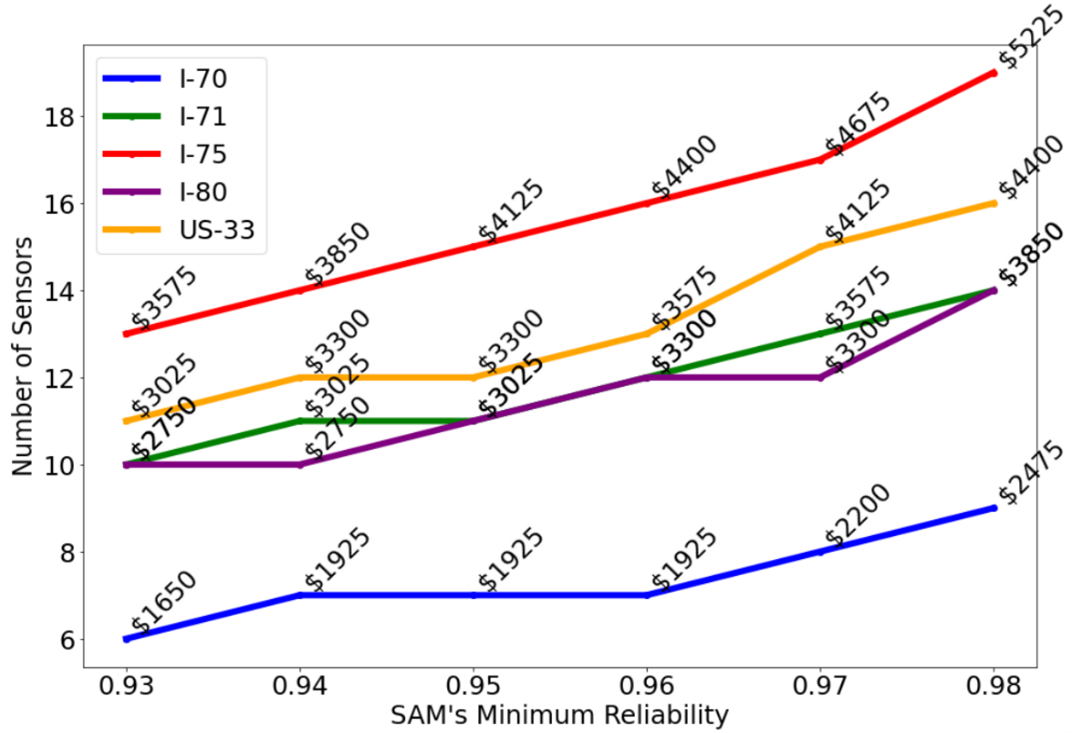}
\caption{Number of ADSB sensors suggested by the model along five corridors under varying SAM’s minimum reliability, with total cost at each point.}
\label{1r_result1}
\end{figure}

\begin{figure}[tb!]
\centering
\includegraphics[width=18cm,height=18cm]{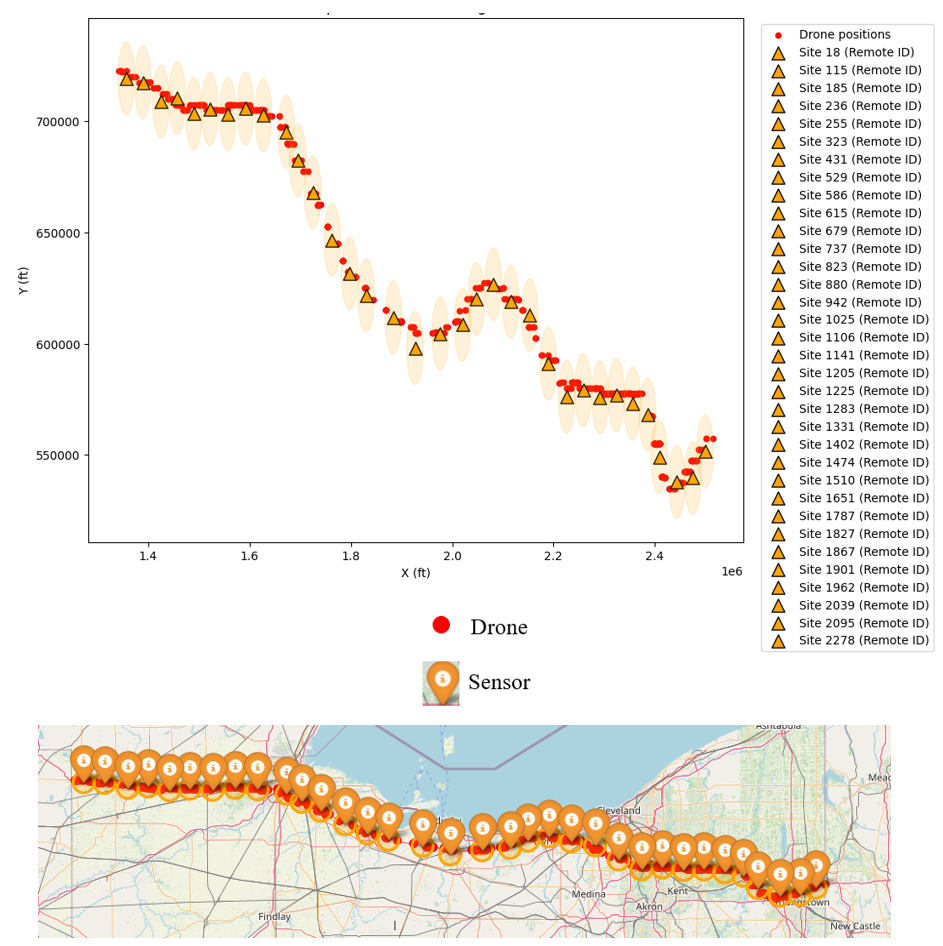}
\caption{Optimal placement of Remote ID sensors along I-80 based on the estimated 2030 flight schedule.}
\label{1r_result2}
\end{figure}

When all six candidate sensor types, ADS-B, Remote ID, Radar, Acoustic, Optical Camera, and RF, are included, the model selects ADS-B as the optimal sensor type across all corridors. Figure \ref{1r_result1} presents the number of ADS-B sensors suggested by the model for the five corridors under varying minimum reliability thresholds, along with the corresponding total cost. As the reliability threshold increases, the total cost of the sensor network also rises because more sensors are required to satisfy the higher reliability target. This trend reflects the fundamental trade-off between safety assurance and infrastructure investment. However, ADS-B may not be a feasible sensor option for AAM operations. AAM vehicles are not required to carry ADS-B transponders. This means they may not broadcast ADS-B signals at all. If the SAM network were built with ground-based ADS-B receivers, those receivers would only be able to detect aircraft that transmit ADS-B. As a result, the system would fail to detect AAM aircraft that do not carry ADS-B equipment. Therefore, even though ADS-B appears optimal from a cost and reliability standpoint in the model, it cannot serve as a reliable sensing solution for monitoring AAM traffic. To account for this operational constraint, we next excluded ADS-B from the candidate set. When the model was rerun with a reliability threshold of 0.95 and ADS-B removed, Remote ID emerged as the optimal sensor type. Figure \ref{1r_result2} shows an example of the resulting optimal Remote ID sensor placement along the I-80 corridor based on the estimated 2030 flight schedule.

\begin{figure}[tb!]
\centering
\includegraphics[width=15cm,height=7.5cm]{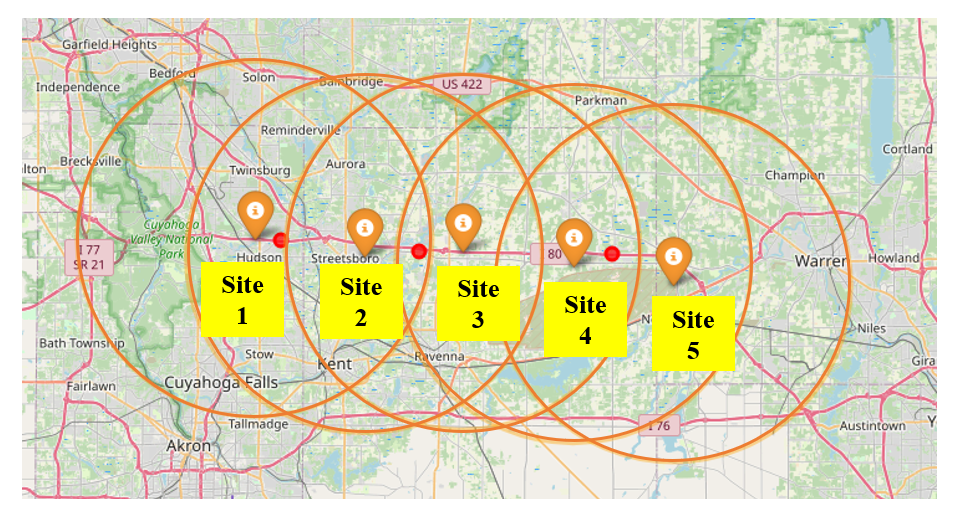}
\caption{Weather-affected region along the I-80 corridor with existing sensors (orange) determined by the reliability model under normal traffic conditions (three aircraft in the schedule).}
\label{2r_results1}
\end{figure}

\begin{figure}[tb!]
\centering
\includegraphics[width=15cm,height=7.5cm]{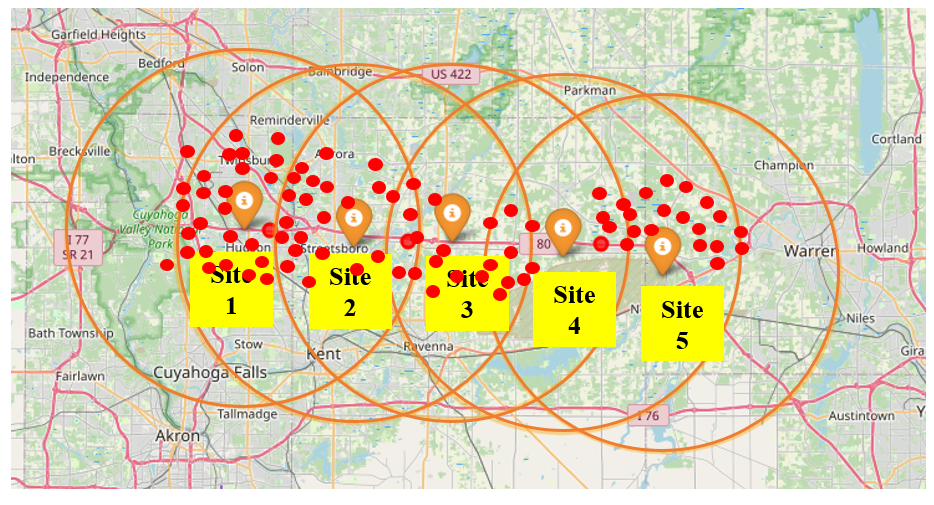}
\caption{Weather-affected region along the I-80 corridor with existing sensors (orange) from the reliability model under increased traffic conditions, where the number of aircraft rises from three to 91 due to congestion.}
\label{2r_results2}
\end{figure}

\begin{figure}[tb!]
\centering
\includegraphics[width=15cm,height=7.5cm]{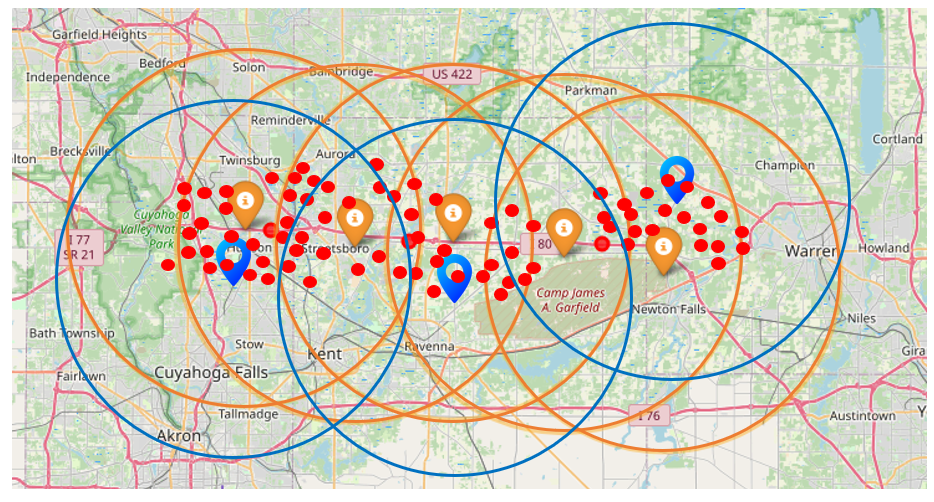}
\caption{Additional sensors (blue) identified by the robustness model to support increased traffic demand, shown alongside the existing sensors (orange).}
\label{2r_results3}
\end{figure}

\subsection{Robustness Model}

The reliability and robustness model work together to guide the design of the SAM sensor network under both normal and perturbed conditions. The reliability model first determines the number, type, and placement of sensors required under normal 2030 AAM traffic conditions, which serves as the baseline or “existing” sensor network. When external perturbations occur, such as sudden increases in AAM traffic, robustness model A identifies the additional sensors needed to maintain a minimum SAM performance. The robustness model A allows to specify the corridor of interest, define the area to analyze using boundary coordinates, and indicate the number of drones to simulate in the scenario. Using these inputs, the model generates a map showing the optimal placement of additional sensors, including the number of sensors at each location and their types.

The reliability and robustness models together support the design of the SAM sensor network under both normal and perturbed conditions. The reliability model determines the number, type, and placement of sensors under normal traffic conditions, which serves as the baseline (“existing”) network. Figure \ref{2r_results1} shows one illustrative example of this baseline in a selected test area, where the orange icons represent the sensors and the red markers indicate aircraft operating under regular traffic levels. In this case, five existing sensors are sufficient to provide surveillance under normal conditions, where approximately three aircraft are present on average at each time step.

To assess system performance under disruptions, we consider scenarios where traffic demand increases sharply due to external factors such as adverse weather, flight cancellations, or rerouting. For example, the number of aircraft increases from 3 to 91, as shown in Figure \ref{2r_results2}. This sudden surge creates congestion that can exceed the monitoring capacity of the existing sensor network. In this scenario, to maintain a minimum level of surveillance performance, the robustness model identifies three additional sensors that need to be activated under these degraded conditions. Using the baseline sensors from the reliability model as input, it determines the optimal locations and types of these sensors. Figure \ref{2r_results3} shows the output, where orange icons represent the existing sensors and blue icons indicate the three additional sensors activated to handle the increased traffic demand.


\subsection{Resiliency Model}

\begin{figure}[tb!]
\centering
\includegraphics[width=14cm,height=9cm]{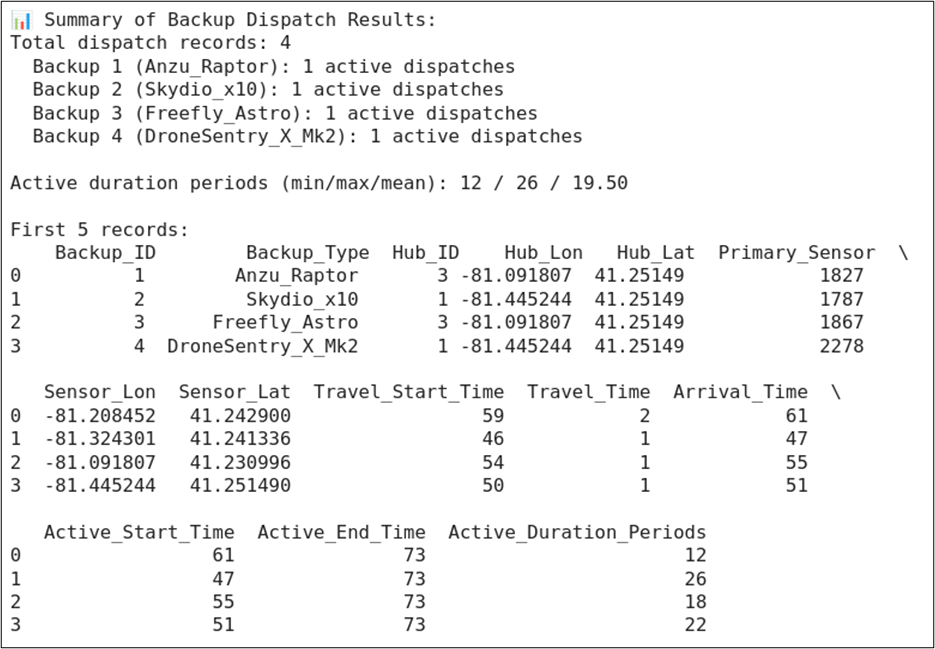}
\caption{Dispatch summary of backup sensors.}
\label{3r_result1}
\end{figure}

\begin{figure}[tb!]
\centering
\includegraphics[width=17cm,height=8cm]{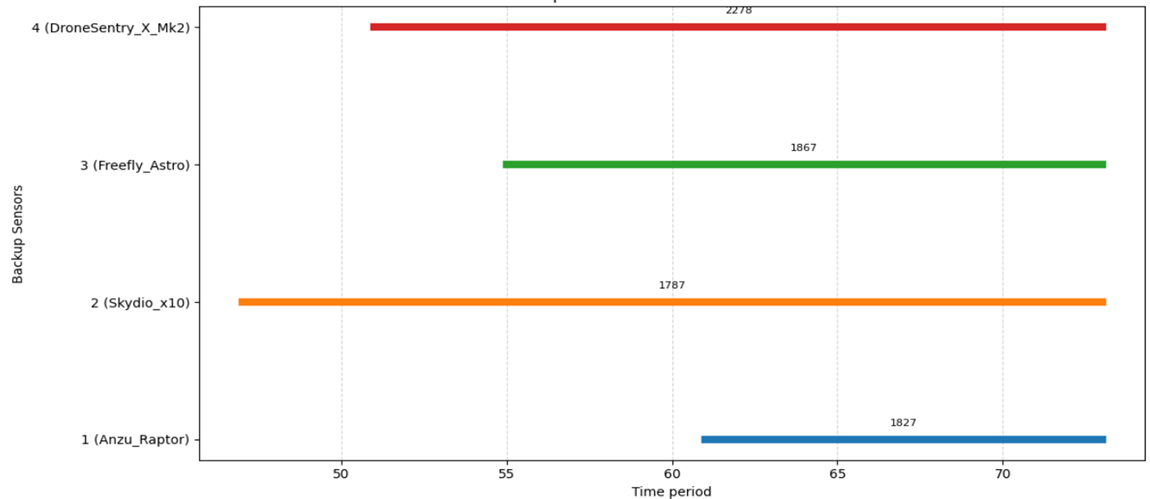}
\caption{Visualization of backup sensor deployment schedule showing when each backup unit is active and the corresponding primary sensor site IDs they support.}
\label{3r_result2}
\end{figure}

The resiliency model is developed to be used after robustness model B, once the redundancy loss in the sensor network has been evaluated and found to remain below an acceptable threshold. In such cases, the user can run the resiliency model  to determine the optimal backup sensor dispatch strategy, ensuring that the SAM network can maintain or restore performance under disruptive events. The model requires to provide inputs such as the backup UAV or ground vehicle type, hub locations, and the sensors they are intended to support. Based on these inputs, the model generates an optimized schedule for the backup units, specifying the departure hub, travel start time, travel duration, arrival time, active period start and end times, and the total number of periods during which each backup sensor is active. For example, in one scenario as shown in Figure \ref{3r_result1}, the model scheduled four backup units: Anzu Raptor, Skydio x10, Freefly Astro, and DroneSentry X Mk2. Each unit was assigned to a specific hub and sensor target, with corresponding travel and activation times to ensure coverage. Anzu Raptor was dispatched from Hub 3 to support sensor at site ID 1827, starting at time 59, traveling for 2 periods, and remaining active from period 61 to 73. Similarly, Skydio x10 departed from Hub 1 to support sensor at site ID 1787, Freefly Astro from Hub 3 to sensor at site ID 1867, and DroneSentr X Mk2 from Hub 1 to sensor at site ID 2278. Figure \ref{3r_result2} visualizes the deployment schedule of backup sensors generated by the resiliency model. The figure shows when each backup unit is active and the primary sensor site IDs they are assigned to support. This allows users to clearly see the timing and coverage of backup resources.

\section{Conclusion and Future Work}

AAM operations are heavily based on continuous situational awareness in low-altitude airspace. To support safe passenger transport, emergency response, and cargo delivery, surveillance systems must detect and track aircraft reliably—even when components degrade, weather conditions vary, or the system experiences unexpected failures. This study addresses this need by developing the foundation of the 3R framework, consisting of reliability, robustness, and resilience models for the SAM system. In this study, we present the reliability model, which determines the optimal placement of multi-type sensors to meet the required detection thresholds under normal operating conditions. Using DSM/LiDAR terrain data, corridor-level AAM demand forecasts, and detailed sensor specifications, the model identifies the minimum-cost deployment that ensures adequate coverage. The robustness model identifies additional sensor requirements under external perturbations, such as increased AAM traffic, and the resilience model develops backup deployment strategies to maintain temporary surveillance during primary sensor outages.

Future work will extend the proposed 3R framework by considering additional failure and uncertainty scenarios that may affect the performance of the SAM system. For example, future studies can include communication-related failures such as data transmission errors, packet loss, network delays, server outages, and unstable communication links to evaluate how these disruptions impact overall surveillance reliability. Sensitivity analysis can also be conducted by varying important input parameters and evaluating how the model performs under different operating conditions. These parameters may include sensor reliability, communication failure rates, traffic demand, sensing range, weather conditions, and reliability thresholds. Such analyses would help assess the stability and adaptability of the proposed framework and identify the most influential factors affecting system performance. Future work can also investigate extreme disruption scenarios in which complete surveillance outages occur and backup sensors are unavailable for deployment. In such situations, maintaining safe AAM operations becomes significantly more challenging because the system may temporarily lose the ability to reliably detect and track aircraft within affected regions. To address this issue, a mitigation and recovery model can be developed to support operational decision-making during large-scale failures. It can determine alternative operational strategies such as temporary corridor closures, dynamic flight rerouting, altitude restrictions, reduced traffic density, emergency surveillance sharing from neighboring regions, or priority-based allocation of limited sensing resources to critical aircraft operations, such as emergency response vehicles. The framework can also incorporate adaptive risk-based control policies that continuously evaluate the remaining surveillance capability and recommend safe operational limits until normal system performance is restored. Such extensions would further strengthen the resilience of the SAM framework under severe and low-probability disruption events.

\bibliography{sample}

\end{document}